\documentclass[12pt]{article}
\usepackage{latexsym,amssymb,amsmath,fancyhdr,hyperref}

\newtheorem{theorem}{Theorem}[section]

\newtheorem{corollary}[theorem]{Corollary}

\newtheorem{lemma}[theorem]{Lemma}

\newtheorem{definition}[theorem]{Definition}
\newtheorem{example}[theorem]{Example}

\newcommand{\R}{\mathbb{R}}

\newcommand{\T}{\mathbb{T}}
\newcommand{\D}{\Delta}
\newcommand{\si}{\sigma}
\newcommand{\cR}{{\cal R}}
\newcommand{\ka}{{\kappa}}

\newcommand{\cF}{\cal F}
\newcommand{\cS}{\cal S}
\newcommand{\cC}{\cal C}

\title{On the non--multiplicity of solutions to \\ matrix equations on time scales \footnote{Submitted to \textit{Electronic Journal of Differential Equations} on November 13, 2012}}

\author{Atiya H. Zaidi \footnote{{\bf Corresponding Author}: a.zaidi@unsw.edu.au}}
\date{}

\begin{document}

\maketitle
\begin{center}
School of Mathematics and Statistics,\\
        The University of New South Wales,\\
         Sydney NSW 2052, Australia 
\end{center}

\begin{abstract}
In this paper we establish the non--multiplicity of solutions to first order matrix dynamic equations on time scales. The new results verify and extend the notions developed in \cite{thesis} to more complex systems of $n^2$ matrices with the help of ideas developed in \cite[Chap 5]{BP}, identifying Lipschitz conditions suitable to generalised $n^2$--models on time scales.
\end{abstract}




\section{Introduction}
\setcounter{equation}{0}
\setcounter{section}{1}

The study of dynamic equations on time scales was initiated in 1988 by S. Hilger when he introduced the concept and the calculus of unifying mathematical analyses of continuous and discrete dynamics, see \cite{Hilger}, \cite{Hilger2}. Since then, several results have been developed to complement his ideas to shape the linear and the nonlinear theory of dynamic equations on time scales. These equations describe continuous, discrete or both types of phenomena occurring simultaneously, through a single model. 

In \cite{TZ1} and \cite{thesis} we presented results regarding non--multiplicity of solutions to nonlinear models of dimension $n$ on time scales. 
In this work we plug in some of those notions to understand more complex systems of dimension $n^2$ for $n \ge 1$. Most physical processes that occur in nature, industry and society are nonlinear in structure and depend on several factors and their interactions. Also, in real life problems, it may not be possible to change the initial or prevailing states of a dynamic model as well as the natural or circumstantial relationships of the variables involved. Knowing that a mathematical formulation of such a system with the given initial conditions has either one solution or no solution would lead to the guarantee that `existence' of a solution implies its uniqueness. 

\indent This paper considers two basic types of dynamic IVPs of dimension $n^2$. These are:

\begin{eqnarray}
X^\D &=& F(t, X);  \label{1}\\
\mbox{and} \qquad \ X^\D &=& F(t, X^\si), \label{2} 
\end{eqnarray}
subject to the initial condition
\begin{equation}
X(a) = A. \label{1i}
\end{equation}

In the above systems, $X$ is a $n^2$--matrix--valued function on a time scale interval $[a,b]_\T:=[a,b] \cap \T$, where $\T$ is a non-empty and closed subset of $\R$, with $b > a$; $F:[a,b]_\T \times \R^{n^2} \rightarrow \R^{n^2}$; $X^\si = (x^\si_{ij})$ and $X^\D = (x^\D_{ij})$ for $1\le i,j \le n$; and $A$ is a given constant $n^2$ matrix. A solution of \eqref{1}, \eqref{1i} (respectively \eqref{2}, \eqref{1i}) will be a matrix--valued function $X$ which solves \eqref{1} and \eqref{1i} (respectively \eqref{2}, \eqref{1i}) on $[a,b]_{\T}$.

Our main aim in this work is to derive conditions that would ensure that there is either one or no solution to initial value problems \eqref{1}, \eqref{1i} and \eqref{2}, \eqref{1i}. Our new results significantly improve those in \cite{thesis} and present some novel ideas.

In the next section, we identify some basic concepts of the time scale calculus associated with matrix--valued functions, used in this work.

\section{Preliminaries}\label{sec2}
\setcounter{equation}{0}
\setcounter{section}{2}

The following definitions and descriptions explain how we use the time scale notation within the set of $m \times n$ matrices on $\T$. For more detail see \cite{BS}, \cite{BP}, \cite{Hilger},  \cite{KRP}, \cite{thesis}.

\begin{definition} 
Let $\T$ be an arbitrary time scale and $t$ be a point in $\T$. The forward jump operator, $\si(t): \T \to \T$, is defined as $\si (t):= \inf \lbrace s \in \T: s > t \rbrace$ for all $t \in \T.$ In a similar way, we define the backward jump operator, $\rho(t): \T \to \T$, as $\rho (t):= \sup \lbrace s \in \T: s < t \rbrace$ for all $t \in \T.$
\end{definition}
\begin{flushright}
$\square$
\end{flushright}

In this way, the forward and backward (or right and left) jump operators declare whether a point in a time scale is discrete and give the direction of discreteness of the point. The results in this paper concern the forward or rightward motion on $[a,b]_\T$. Hence, further notation and definitions will be presented accordingly.

Continuity of a function at a point $t \in \T$ is said to be `right--dense' when $t = \si(t)$, otherwise it is called right--scattered. The `step size' at each point of a time scale is given by the graininess function, $\mu(t)$, defined as $\mu(t):= \si(t) - t$ for all $t \in \T$. If $\T$ is discrete, it has a left--scattered maximum value $m$ and we define $\T^{{\kappa}}:= \T \setminus m$, otherwise $\T^{{\kappa}} := \T$. 

%
%

Analogous to left--Hilger--continuous functions \cite[p.3]{zaidi} for any ordered $n$--pair $(t, {\bf x}) \in \T \times \R^n$, we define a \emph{right--Hilger--continuous} function ${\bf f}(t, {\bf x})$ \cite{Hilger}, \cite[Chap.2]{thesis} as a function ${\bf f}: \T^\kappa \times \R^{n} \to \R^n$ having the property that ${\bf f}$ is continuous at each $(t,{\bf x})$ where $t$ is right--dense; and the limits $$ \lim_{(s,{\bf y}) \to (t^{-},{\bf x})} {\bf f}(s, {\bf y}) \qquad \mbox{and} \qquad \lim_{{\bf y} \to {\bf x}} {\bf f}(t,{\bf y})$$ both exist and are finite at each $(t,{\bf x})$ where $t$ is left--dense.

It should be noted that ${\bf f}$ is \emph{rd--continuous} if ${\bf f}(t, {\bf x}) = {\bf g}(t)$ for all $t \in \T$ and is \emph{continuous} if ${\bf f}(t,{\bf x}) = {\bf h}({\bf x})$ for all $t \in \T$.\\

Continuity of a matrix--valued function at a point $t \in \T$ depends on the continuity of its elements at $t$. Thus, for any $t \in \T$, a rd--continuous matrix--valued function is a function $X: \T \rightarrow \R^{m \times n}$ with entries $(x_{ij})$, where $x_{ij}:\T \to \R$; $1\le i \le m, 1 \le j \le n$; and each $x_{ij}$ is rd--continuous on $\T$. Moreover, we say that $X \in C_{rd}=C_{rd}(\T;\R^{m \times n})$ \cite[p.189]{BP}.\\

Thus, a right--Hilger--continuous matrix--valued function can be defined as follows.

\begin{definition} \label{rhcont} 
Assume $F: \T \times \R^{m \times n} \to \R^{m \times n}$ be a matrix--valued function with entries $(f_{ij})$, where each $f_{ij}:\T \times \R \to \R$ for $1\le i \le m, 1 \le j \le n$. We define $F$ to be right--Hilger--continuous if each $f_{ij}(t,x_{kl})$ is right--Hilger--continuous for all $t \in \T$ and $x_{kl}:\T \to \R$ for all $k,l$.
\end{definition}
\begin{flushright}
$\square$
\end{flushright}


For a fixed $t \in \T^{{\kappa}}$ and $x: \T \rightarrow \R$, the delta--derivative of $x$ (if it exists) is $x^\D(t)$, having the property that given $\epsilon > 0$ there is a neighbourhood $U$ of $t$, that is, $U = (t - \delta, t + \delta) \cap \T$ for some $\delta > 0$, such that $$ |(x^\si(t) - x(s)) - x^\D (\si(t) - s)| \le \epsilon |\si(t) - s|, \mbox{for all} \ s \in U.$$ Hence, the delta--derivative of a matrix--valued function on a time scale is defined as follows.

\begin{definition} 
Consider a function $X:\T \to \R^{m \times n}$. We define $X^{\D}:= (x_{ij}^\D)$  to be the delta--derivative of $X$ on $\T$ if $x_{ij}^\D(t)$ exists for all $t \in \T^{\kappa}$ for all $1 \le i \le m, 1 \le j \le n$ and say that $X$ is delta--differentiable on $\T$.
\end{definition}
\begin{flushright}
$\square$
\end{flushright}

The set of delta--differentiable matrix--valued functions $K: \T \rightarrow \R^{m \times n}$ satisfy the simple useful formula \cite[Theorem 5.2]{BP}
\begin{equation}
K^{\si}(t) = K(t) + \mu(t)K^{\D}(t),  \qquad \mbox{for all} \ t \in \T^{\kappa}. \label{suf}
\end{equation}

The next theorem describes some more identities related to delta--differentiable matrix--valued functions that will be used in this work \cite[Theorem 5.3]{BP}.

\begin{theorem} \label{main} 
Let $X,Y: \T \rightarrow \R^{n^2}$ be matrix--valued functions. If $X, Y$ are delta--differentiable on $\T$ then for all $t \in \T^{\kappa}$ we have  
\begin{enumerate}
\item[a.] $(X+Y)^\D(t) = X^\D(t) + Y^\D(t)$;
\item[b.] for any constant $k \in \R$, $(kX)^\D(t) = k X^\D(t)$;
\item[c.] $(XY)^\D(t) = [X^\D Y + X^\si Y^\D](t) = [X Y^\D + X^\D Y^\si](t)$;
\item[d.] If $X(t)$ and $X^\si(t)$ are invertible for all $t \in \T^{\kappa}$ then $$(X^{-1})^\D(t) = [-X^{-1} X^\D (X^\si)^{-1}](t) = [-(X^\si)^{-1} X^\D X^{-1}](t);$$
\item[e.] If $Y(t)$ and $Y^\si(t)$ are invertible for all $t \in \T^{\kappa}$ then $$[XY^{-1}]^\D(t) = [X^\D - X Y^{-1} Y^\D](t) (Y^\si(t))^{-1} = [X^\D - (XY^{-1})^\si Y^\D](t) Y^{-1}(t);$$
\item[f.] $(X^{*})^{\D} = (X^{\D})^{*}$, where $*$ refers to the conjugate transpose.  
\end{enumerate}
\end{theorem}
\begin{flushright}
$\square$
\end{flushright}


Since all rd--continuous functions are delta--integrable, the antiderivative of a right--Hilger--continuous matrix--valued function can be defined as follows:

\begin{theorem} \label{intf1} 
Let $F: \T^{\kappa} \times \R^{n^2} \to \R^{n^2}$ and $a \in \T$. If $F$ is right--Hilger--continuous on $\T^{\kappa} \times \R^{n^2}$ then there exists a function ${\cF}:C(\T;\R^{n^2}) \to C(\T;\R^{n^2})$ called the delta integral of $F$ such that
\begin{eqnarray} \label{Hintegral}
[{\cF}X](t):= \int_{a}^{t} F(s, X(s)) \ \D s, \qquad \mbox{for all} \ t \in \T.
\end{eqnarray}
\end{theorem}
\begin{flushright}
$\square$
\end{flushright}

Next, we describe positive definite (respectively semi--definite) $n^2$--matrices and some of their properties \cite{bhatia}, \cite{Horn}, \cite{KRP} on a time scale $\T$. This class of square matrices on $\T$ plays a vital role in establishing the non--multiplicity of solutions in this work.
 
\begin{definition} \label{posdef} 
Let $X:[a,b]_\T \to \R^{n^2}$ and ${\bf z}:\T \to \R^n$. Assume $\bf z \neq 0$ for all $t \in [a,b]_\T$. We say that $X$ is positive definite (respectively semi--definite) on $[a,b]_\T$ if ${\bf z}^TX{\bf z} > 0$ (respectively ${\bf z}^TX{\bf z} \ge 0$) on $[a,b]_\T$ and write $X > 0$ (respectively $X \ge 0$) on $[a,b]_\T$.
\end{definition}
\begin{flushright}
$\square$
\end{flushright}

It is clear from the above definition that a negative definite (respectively semi--definite) matrix $Y$ on $\T$ will satisfy ${\bf z}^{T}Y{\bf z} < 0$ (respectively ${\bf z}^{T}Y{\bf z} \le 0$) for all ${\bf z}:\T \to \R^n$ and we say that $Y < 0$ (respectively $Y \le 0$). 

The class of positive definite matrices defined above has the following properties. 

\begin{theorem} \label{propposdef}
Let $A, B:[a,b]_\T \to \R^{n^2}$. If $A, B > 0$ on $[a,b]_\T$ then the following properties hold on $[a,b]_\T$:
\begin{enumerate}
\item $A$ is invertible and $A^{-1} > 0$;
\item if $\alpha \in \R$ such that $\alpha > 0$ then $\alpha A > 0$;
\item if $\lambda$ is an eigenvalue of $A$ then $\lambda > 0$;
\item $det(A) > 0$ and $tr(A)>0$.
\item $A+B > 0$, $ABA > 0$ and $BAB > 0$;
\item if $A$ and $B$ commute then $AB > 0$ and similarly, if there exists ${\cC} \le 0$ such that $A$ and $\cC$ commute then $AC \le 0$;
\item if $A - B \ge 0$ then $A \ge B$ and $B^{-1} \ge A^{-1} > 0$;
\item there exists $\beta > 0$ such that $A > \beta I$.
\end{enumerate} 
\end{theorem}
\begin{flushright}
$\square$
\end{flushright}

The regressiveness of $n^2$--matrix functions and their properties is described \cite{BP} in a similar manner as for regressive $n$--functions, as follows.
 
\begin{definition} 
Consider a function $K:\T \to \R^{n^2}$. We call $K$ regressive on $\T$ if the following conditions hold:
\begin{itemize}
\item $K$ is rd--continuous on $\T$; and
\item the matrix $I+\mu(t)K$ is invertible for all $t \in \T^{\kappa}$, where $I$ is the identity matrix.
\end{itemize}
We denote by $ \cR:= \cR(\T;\R^{n^2})$ the set of all regressive $n^2$--matrix functions on $\T$.
\end{definition}
\begin{flushright}
$\square$
\end{flushright}

It is clear from above that all positive and negative definite matrix--valued functions are regressive.  The following theorem \cite[pp. 191-192]{BP} lists some important properties of regressive $n^2$--matrix functions on $\T$.

\begin{theorem} \label{propoplus}
Let $A,B:\T \to \R^{n^2}$. If $A,B \in \cR$ then the following identities hold for all $t \in \T^{\kappa}$:
\begin{enumerate}
\item $(A \oplus B)(t) = A(t) + B(t) + \mu(t)A(t)B(t)$;
\item $(\ominus A)(t) = -[I+\mu(t)A(t)]^{-1} A(t) = -A(t)[I+\mu(t)A(t)]^{-1}$;
\item $A^{*} \in \cR$ and $(\ominus A)^{*}=\ominus A^{*}$;
\item $I+\mu(t)(A(t)\oplus B(t)) = [I+\mu(t)A(t)][I+\mu(t)B(t)]$;
\item $I+\mu(t)(\ominus A(t))= [I+\mu(t)A(t)]^{-1}$;
\item $(A \ominus B)(t) = (A \oplus (\ominus B))(t) = A(t)-[I+\mu(t)A(t)][I+\mu(t)B(t)]^{-1} B(t)$;
\item $[A(t) \oplus B(t)]^{*} = A(t)^{*} \oplus B(t)^{*}$.
\end{enumerate} 
\end{theorem}
\begin{flushright}
$\square$
\end{flushright}

An important implication of regressive matrices is the generalised matrix exponential function on a time scale.

\begin{definition} \label{exp}
Let $K:\T \to \R^{n^2}$ be a matrix--valued function. Fix $a \in \T$ and assume $P \in \cR$. The matrix exponential function denoted by $e_{K}(\cdot, a)$ is defined as
\begin{eqnarray}
e_{K}(t,a):= \left\{
\begin{array}{lr}
\exp \left(\int_{a}^{t} K(s) \ d s \right), & \mbox{for} \ t \in \T, \ \mu = 0; \\
\\
\exp \left(\int_{a}^{t} \dfrac{\mbox{Log} (I + \mu(s) K(s))}{\mu(s)} \ \D s \right), & \mbox{for} \ t \in \T, \mu > 0,
\end{array} \right.
\end{eqnarray}
where $\mbox{Log}$ is the principal logarithm function.
\end{definition}
\begin{flushright}
$\square$
\end{flushright}

Further properties of the matrix exponential function \cite[Chap 5]{BP} are shown in the following theorem and will be used in this work.

%

\begin{theorem} \label{propexp} 
Let $K,L:\T \to \R^{n^2}$. If $K,L \in \cR$ then the following properties hold for all $t,s,r \in \T$:
\begin{enumerate}
\item $e_{0}(t,s) = I = e_{K}(t,t)$, where $0$ is the $n^2$ zero matrix;
\item $e_{K}^{\si}(t,s) = e_{K}(\si(t),s) = (I + \mu(t)K(t))e_{K}(t,s)$;
\item $e_{K}(s,t) = e_{K}^{-1}(t,s) = e_{\ominus K^*}^{*}(t,s)$;
\item $e_{K}(t,s) e_{K}(s,r) = e_{K}(t,r)$;
\item $e_{K}(t,s)e_{L}(t,s) = e_{K \oplus L}(t,s)$;
\item $e_{K}^{\D}(t,s) = -e_{K}^{\si}(t,s) K(t) = K(t) e_{K}(t,s)$.
\end{enumerate}
\end{theorem}
\begin{flushright}
$\square$
\end{flushright}

\section{Lipschitz continuity of matrix functions on $\T$} \label{sec3}
\setcounter{equation}{0}
\setcounter{section}{3}

In this section, we present Lipschitz conditions for matrix--valued functions defined on a subset of $\T \times \R^{n^2}$ that allow positive definite matrices as Lipschitz constants for these functions. The ideas are obtained from \cite{AL}, \cite{EC} \cite{Hart} \cite{KP} and \cite{thesis}.

\begin{definition} \label{leftsided} 
Let $S \subset \R^{n^2}$ and $F: [a,b]_{\T} \times S \to \R^{n^2}$ be a right--Hilger--continuous function. If there exists a positive definite matrix $B$ on $\T$ such that for all $P, Q \in S$ with $P > Q$, the inequality 
\begin{equation} \label{LCOS}
F(t,P) - F(t,Q) \le B(t) (P-Q), \qquad \mbox{for all} \ (t,P), (t,Q) \in [a,b]]^{\kappa}_{\T} \times S
\end{equation}
holds, then we say $F$ satisfies a left--handed--Lipschitz condition (or is left--handed Lipschitz continuous) on $[a,b]_{\T} \times S$. 
\end{definition}
\begin{flushright}
$\square$
\end{flushright}

\begin{definition} \label{rightsided} 
Let $S \subset \R^{n^2}$ and $F: [a,b]_{\T} \times S \to \R^{n^2}$ be a right--Hilger--continuous function. If there exists a positive definite matrix $\cC$ on $\T$ such that for all $P, Q \in S$ with $P > Q$, the inequality 
\begin{equation} \label{LCOSright}
F(t,P) - F(t,Q) \le (P-Q) C(t), \qquad \mbox{for all} \ (t,P), (t,Q) \in [a,b]]^{\kappa}_{\T} \times S
\end{equation}
holds, then we say $F$ satisfies a right--handed--Lipschitz condition (or is right--handed Lipschitz continuous) on $[a,b]_{\T} \times S$. 
\end{definition}
\begin{flushright}
$\square$
\end{flushright}

Classically, any value of matrix $B$ or $\cC$ satisfying \eqref{LCOS} or \eqref{LCOSright} would depend only on $[a,b]_{\T} \times S$ \cite[p.6]{Hart}. For the sake of simplicity, we consider $[a,b]^{\kappa}_\T \times S$ to be convex and $F$ smooth on $[a,b]^{\kappa}_\T \times S$, then the following theorem \cite[p.248]{EC}, \cite[Lemma 3.2.1]{AL} will be helpful to identify a Lipschitz constant for $F$ on $[a,b]^{\kappa}_\T \times S$ and obtain a sufficient condition for $F$ to satisfy the left-- or right--handed Lipschitz condition on $[a,b]^{\kappa}_\T \times S$.

\begin{corollary} \label{exisLC}
Let $a, b \in \T$ with $b > a$ and $A \in \R^{n^2}$. Let $k>0$ be a real constant and consider a function $F$ defined either on a rectangle
\begin{equation}
R^{\kappa}:= \{(t,P) \in [a,b]^{\kappa}_{\T} \times \R^{n^2}: \|P - A\| \le k \}
\end{equation}
or on an infinite strip 
\begin{equation}
{\cS^\kappa}:= \{(t,P) \in [a,b]^{{\kappa}}_{\T} \times \R^{n^2}: \|P\| \le \infty\} 
\end{equation}
If $\dfrac{\partial F(t, P)}{\partial p_{ij}}$ exists for all $1 \le i,j \le n$ and is continuous on $R^{{\kappa}}$ (or $\cS^{\kappa}$), and there is a positive definite matrix $L$ such that for all $(t, P) \in R^{{\kappa}}$ (or $\cS^{\kappa}$), we have
\begin{equation} \label{parder}
\dfrac{\partial F(t, P)}{\partial p_{ij}} \le L, \ \qquad \mbox{for all} \ i,j = 1,2,\cdots,
\end{equation}
then $F$ satisfies \eqref{LCOS} with $B(t) = L$ or \eqref{LCOSright} with ${\cC}(t) = L$, on $R^{{\kappa}}$ (or $\cS^{\kappa}$) for all $t \in [a,b]_{\T}$.
\end{corollary}

\noindent {\bf Proof}: The proof is similar to that of \cite[Lemma 3.2.1]{AL} except that $\dfrac{\partial F(t, P)}{\partial p_{ij}}$ is considered bounded above by $B(t) = L$ in the left--handed case or ${\cC}(t) = L$ in the right--handed case, for all $t \in [a,b]_{\T}$.
\begin{flushright}
$\square$
\end{flushright}

\section{non--multiplicity results}
\setcounter{equation}{0}
\setcounter{section}{4}

In this section, we present generalised results regarding non--multiplicity of solutions to the dynamic IVPs \eqref{1}, \eqref{1i} and \eqref{2}, \eqref{1i} within a domain $S \subseteq \R^{n^2}$. The results are heavily based on ideas in \cite[Chap 5]{BP}, methods from ordinary differential equations \cite{EC}, \cite{Birk} and \cite{KP} and matrix theory \cite{bhatia}, \cite{Horn}, \cite{pipes}.

The following lemma establishes a function to be a solution of \eqref{1}, \eqref{1i} and \eqref{2}, \eqref{1i}.

\begin{lemma} \label{soldeqX}
Consider the dynamic IVP \eqref{1}, \eqref{1i}. Let $F: [a,b]^{\kappa}_{\T} \times \R^{n^2} \rightarrow \R^{n^2}$ be a right--Hilger--continuous matrix--valued function. Then a function $X$ solves \eqref{1}, \eqref{1i} if and only if it satisfies
\begin{eqnarray} \label{deqsolX}
X(t) = \int_{a}^{t} F(s, X(s)) \ \D s + A, \qquad \qquad \qquad \mbox{for all} \ t \in [a,b]_{\T}.
\end{eqnarray}
\end{lemma}
\begin{flushright}
$\square$
\end{flushright}

A similar function can be defined as a solution of \eqref{2},\eqref{1i}. \\

\begin{theorem} \label{exis1} 
Let $S \subseteq \R^{n^2}$ and let $F:[a,b]_\T \times S \rightarrow \R^{n^2}$ be a right--Hilger--continuous function. If there exist $P, Q \in S$ with $P > Q$ and a positive definite matrix $B$ on ${\T}$ such that
\begin{enumerate}
\item[(1)] $B \in C_{rd}([a,b]_{\T};\R^{n^2})$;
\item[(2)] $e_{B}(t,a)$ commutes with $B(t)$ for all $t \in [a,b]_{\T}$ and with $P(t)$ for all $(t,P) \in [a,b]_{\T} \times S$;
\item[(3)] the left--handed Lipschitz condition, $F(t,P) - F(t,Q) \le B(t) (P-Q)$ holds for all $(t,P), (t,Q) \in [a,b]]^{\kappa}_\T \times S$,
\end{enumerate}
then the IVP \eqref{1}, \eqref{1i} has, at most, one solution, $X$, with $X(t) \in S$ for all $t \in [a,b]_{\T}$.
\end{theorem}

\noindent {Proof}: We present the proof by contradiction and, without loss of generality, assume two solutions $X, Y$ of \eqref{1}, \eqref{1i} in $S$ such that $X-Y > 0$ on $[a,b]_{\T}$ and show that $X \equiv Y$ on $[a,b]_{\T}$.

By Lemma \ref{soldeqX}, $X$ and $Y$ must satisfy \eqref{deqsolX}. Define $U:= X-Y$ on $[a,b]_{\T}$. We show that $U \equiv 0$ on $[a,b]_{\T}$.

Since $\textit{(3)}$ holds, we have, for all $t \in [a,b]_{\T}^{\kappa}$ 
\begin{equation} \label{LCOS2}
U^{\D}(t) - B(t) U(t) =  F(t, X(t)) - F(t, Y(t)) - B(t) (X(t) - Y(t)) \le 0.
\end{equation}
Note that $B$ being positive definite is regressive on $[a,b]_{\T}$. Thus, we have $e_{B}(t,a), e_{B}^{\si}(t,a)$ positive definite with positive definite inverses on $[a,b]_{\T}$, by Theorem \ref{propposdef}(1). Hence, using Theorem \ref{main} and Theorem \ref{propexp} we obtain, for all $t \in [a,b]_{\T}^{\ka}$,
\begin{eqnarray*}
[e_{B}^{-1}(t,a) U(t)]^{\D} &=& [e_{B}^{-1}(t,a)]^{\si} U^{\D}(t) + [e_{B}^{-1}(t,a)]^{\D} U(t)\\ 
&=& [e_{B}^{\si}(t,a)]^{-1} U^{\D}(t) - [e_{B}^{\si}(t,a)]^{-1} e_{B}^{\D}(t,a) e_{B}^{-1}(t,a) U(t)\\
&=& [e_{B}^{\si}(t,a)]^{-1} [U^{\D}(t) - e_{B}^{\D}(t,a) e_{B}^{-1}(t,a) U(t)]\\
&=& (e_{B}^{\si}(t,a))^{-1} [U^{\D}(t) - B(t) U(t)].
\end{eqnarray*}
By virtue of $\emph{(2)}$, $e_{B}^{-1}(t,a)$ also commutes with $B(t)$ for all $t \in [a,b]_{\T}$ and with $P^{\D}(t)$ for all $(t,P) \in [a,b]_{\T} \times S$. Thus, $e_{B}^{-1}(t,a)$ commutes with $U^{\D}(t) - B(t)U(t)$ for all $t \in [a,b]_{\T}$. Hence, by Theorem \ref{propposdef}(6) and \eqref{LCOS2}, we obtain $$ [e_{B}^{-1}(t,a) U(t)]^{\D} \le 0, \qquad \mbox{for all} \ t \in [a,b]_{\T}^{\kappa}.$$ This means that $e_{B}^{-1}(t,a) U(t)$ is non--increasing for all $t \in [a,b]_{\T}$. But $U$ is positive semi--definite on $[a,b]_{\T}$ and $U(a)=0$. Hence, $U \equiv 0$ on $[a,b]_{\T}$. This means that $X(t) = Y(t)$ for all $t \in [a,b]_{\T}$.
\begin{flushright}
$\square$
\end{flushright}

A similar argument holds for the case where $Y - X \ge 0$ on $[a,b]_\T$.\\

\begin{corollary} \label{exis2} 
The above theorem also holds if $F$ has continuous partial derivatives with respect to the second argument and there exists a positive definite matrix $L$ such that $\dfrac{\partial F(t,P)}{\partial p_{ij}} \le L$. In that case, $F$ satisfies \eqref{LCOS} on $R^{\kappa}$ or $\cS^{\kappa}$ with $B:= L$ by Corollary \ref{exisLC}.
\end{corollary}
\begin{flushright}
$\square$
\end{flushright}

\begin{theorem} \label{exis-right} 
Let $S \subseteq \R^{n^2}$ and let $F:[a,b]_\T \times S \rightarrow \R^{n^2}$ be a right--Hilger--continuous function. If there exist $P, Q \in S$ with $P > Q$ and a positive definite matrix $\cC$ on ${\T}$ such that
\begin{enumerate}
\item[(1)] ${\cC} \in C_{rd}([a,b]_{\T};\R^{n^2})$;
\item[(2)] $e_{\cC}^{-1}(t,a)$ commutes with ${\cC}(t)$ for all $t \in [a,b]_{\T}$ and with $P(t)$ for all $(t,P) \in [a,b]_{\T} \times S$;
\item[(3)] the right--handed Lipschitz condition, $F(t,P) - F(t,Q) \le (P-Q) {\cC}(t)$ holds for all $(t,P), (t,Q) \in [a,b]]^{\kappa}_\T \times S$,
\end{enumerate}
then the IVP \eqref{1}, \eqref{1i} has, at most, one solution, $X$, with $X(t) \in S$ for all $t \in [a,b]_{\T}$.
\end{theorem}

\noindent {Proof}: The proof is similar to Theorem \ref{exis1} and is, therefore, omitted. 
\begin{flushright}
$\square$
\end{flushright}

\begin{corollary} \label{exisright} 
The above theorem also holds if $F$ has continuous partial derivatives with respect to the second argument and there exists a positive definite matrix $H$ such that $\dfrac{\partial F(t,P)}{\partial p_{ij}} \le H$. In that case, $F$ satisfies \eqref{LCOSright} on $R^{\kappa}$ or $\cS^{\kappa}$ with ${\cC}:= H$ by Corollary \ref{exisLC}.
\end{corollary}
\begin{flushright}
$\square$
\end{flushright}

Our next two results are based on the, so called, \textit{inverse Lipschitz condition}, in conjunction with \eqref{LCOS} and \eqref{LCOSright} and determine the existence of at most one solution for \eqref{1}, \eqref{1i} in the light of Theorem \ref{propposdef}(7).

\begin{corollary} \label{exis4} 
Let $S \subseteq \R^{n^2}$ and $F: [a,b]^{\kappa}_{\T} \times S \rightarrow \R^{n^2}$ be right--Hilger--continuous. Assume there exists a positive definite matrix $B$ on $\T$ such that conditions $\emph{(1)}$ and $\emph{(2)}$ of Theorem \ref{exis1} hold. If $P(t)-Q(t)$ is positive definite and increasing for all $(t,P),(t,Q) \in [a,b]_{\T} \times S$ and the inequality
\begin{equation} \label{LCOS4}
(P-Q)^{-1} \le (P^{\D} - Q^{\D})^{-1} B(t) , \qquad \mbox{for all} \ (t,P), (t,Q) \in [a,b]]^{\kappa}_{\T} \times S
\end{equation} 
holds, then the IVP \eqref{1}, \eqref{1i} has, at most, one solution $X$ with $X(t) \in S$ for all $t \in [a,b]_{\T}$.
\end{corollary}

\noindent {\bf Proof:} If \eqref{LCOS4} holds then \eqref{LCOS} holds, by Theorem \ref{propposdef}(7). Hence, the IVP \eqref{1}, \eqref{1i} has, at most, one solution by Theorem \ref{exis1}.
\begin{flushright}
$\square$
\end{flushright}

\begin{corollary} \label{exis4right} 
Let $S \subseteq \R^{n^2}$ and $F: [a,b]^{\kappa}_{\T} \times S \rightarrow \R^{n^2}$ be right--Hilger--continuous. Assume there exists a positive definite matrix ${\cC}$ on $\T$ such that conditions $\emph{(1)}$ and $\emph{(2)}$ of Theorem \ref{exis-right} hold.  If $P(t)-Q(t)$ is positive definite and increasing for all $(t,P), (t,Q) \in [a,b]_{\T} \times S$ and the inequality
\begin{equation} \label{LCOS4right}
(P-Q)^{-1} \le {\cC}(t) (P^{\D} - Q^{\D})^{-1} , \qquad \mbox{for all} \ (t,P), (t,Q) \in [a,b]]^{\kappa}_{\T} \times S
\end{equation} 
holds, then the IVP \eqref{1}, \eqref{1i} has, at most, one solution $x$ with $x(t) \in S$ for all $t \in [a,b]_{\T}$.
\end{corollary}

\noindent {\bf Proof:} If \eqref{LCOS4right} holds then \eqref{LCOSright} holds, by Theorem \ref{propposdef}(7). Hence, the IVP \eqref{1}, \eqref{1i} has, at most, one solution by Theorem \ref{exis-right}.
\begin{flushright}
$\square$
\end{flushright}

We now present examples to reinforce our results proved above.

\begin{example}
Let $S:= \{P \in \R^{2^2}: tr(P^{T}P) \le 2\}$, where $P = \left( \begin{array}{cc} p_{1} & -p_{2} \\ -p_{2} & p_{1} \end{array} \right)$.

Consider the IVP
\begin{eqnarray} \label{ex1}
X^\D = F(t, X) &=& \left( \begin{array}{cc} 1+x_{1}^{2} & t^{2}-x_{2} \\ x_{2} + t & t-x_{1} \end{array} \right), \qquad \mbox{for all} \ t \in [0, b]^{\kappa}_{\T};\\
X(0) &=& I.
\end{eqnarray}
We claim that this dynamic IVP has, at most, one solution, $X$, such that $tr(X^{T}X)\le 2$ for all $t \in [0, b]_{\T}$.
\end{example}

\noindent {\bf Proof}: We show that $F(t,P)$ satisfies the conditions of Theorem \ref{exis1} for all $(t,P) \in [0, b]^{\kappa}_{\T} \times S$. 

Note that for $P \in S$, we have $\sum_{j=1}^{2} p_{j}^{2} \le 1$. Thus, $|p_{j}| \le 1$ for $j=1,2$. Let $L:= \left( \begin{array}{cc} k & 0 \\ 0 & k \end{array} \right)$, where $k \ge 2$, and let ${\bf z} \in \R^2$ such that ${\bf z} = \left( \begin{array}{c} x \\ y \end{array} \right)$, where $x \neq 0 \neq y$. Then, we have
\begin{equation} \label{Lposdef}
{\bf z}^{T} L {\bf z} = k(x^{2} + y^{2}) > 0. 
\end{equation}
Hence, $L$ is positive definite. We note that $F$ is right--Hilger--continuous on $[0,1]^{\kappa}_\T \times S$ as all of its components are rd--continuous on $[0, b]_\T$. Moreover, since $L$ is a diagonal matrix, it commutes with $e_{L}(t,a)$ for all $t \in [0,b]_{\T}$. It can be easily verified that $e_{L}(t,a)$ also commutes with $P$ for all $(t,P) \in [0,b]_{\T} \times S$. 

We show that $F$ satisfies \eqref{LCOS} on $[0,b]^{\kappa}_\T \times S$. Note that for all $t \in [0, b]^{\kappa}_{\T}$ and $P \in S$, we have 
\begin{equation*}
\begin{array}{cc} 
\dfrac{\partial F}{\partial p_{1}} = \left( \begin{array}{cc} 2p_{1} & 0 \\ 0 & -1 \end{array} \right) \quad \mbox{and} \ & \dfrac{\partial F}{\partial p_{2}} = \left( \begin{array}{cc} 0 & -1 \\ 1  & 0 \end{array} \right).
\end{array}
\end{equation*}

Then, we have
\begin{eqnarray*}
{\bf z}^{T}\left(L - \dfrac{\partial F}{\partial p_{1}} \right) {\bf z} &=& (k-2p_{1})x^{2}+(k+1)y^{2}\\
&\ge& (k-2)x^{2}+(k+1)y^{2},
\end{eqnarray*}
and
\begin{eqnarray*}
{\bf z}^{T} \left( L - \dfrac{\partial F}{\partial p_{2}}\right) {\bf z} = k(x^{2} + y^2).
\end{eqnarray*}
Therefore, $L - \dfrac{\partial F}{\partial p_{j}} > 0$ for $j =1, 2$. Hence, by Theorem \ref{propposdef}(7), $\dfrac{\partial F}{\partial p_{j}} < L$ for $j =1, 2$ and employing Corollary \ref{exis2}, we have \eqref{LCOS} holding for $L = \left( \begin{array}{cc} k & 0 \\
0 & k \end{array} \right)$ for all $k \ge 0$. 

In this way, all conditions of Theorem \ref{exis1} are satisfied and we conclude that our example has, at most, one solution, $X(t) \in S$, for all $t \in [0, b]_{\T}$.
\begin{flushright}
$\square$
\end{flushright}

\begin{example}
Let $u,w$ be differentiable functions on $(0, \infty)_{\T}$ with $u$ increasing and $u(t)>1$ for all $t \in (0, \infty)_{\T}$. Let $D$ be the set of all $2^2$-positive definite symmetric matrices. We show that, for any matrix $P$ of the form $P := \left( \begin{array}{cc} 2u+t^{2} & w-t \\ w-t & 2u+t^{2} \end{array} \right)$ in $D$, there exists $Q := \left( \begin{array}{cc} u+t^{2} & w-t \\ w-t & u+t^{2} \end{array} \right)$ in $D$ such that the dynamic IVP \eqref{1}, \eqref{1i} has, at most, one solution, $X$, on $(0, \infty)_{\T}$ such that $X \in D$.
\end{example}

\noindent {\bf Proof}: We show that \eqref{1} satisfies the conditions of Corollary \ref{exis4} for all $(t,P), (t,Q) \in (0, \infty)^{\kappa}_{\T} \times D$. 

Note that since $u,w$ are differentiable on $(0, \infty)_{\T}$, we have $P^{\D} = F(t,P)$ and $Q^{\D} = F(t,Q)$ right--Hilger--continuous on $(0, \infty)^{\kappa}_{\T}$. We also note that $P-Q = \left( \begin{array}{cc} u & 0 \\ 0 & u \end{array} \right)$, which is positive definite and, hence, invertible by Theorem \ref{propposdef}(1). Moreover, since $u^{\D},v^{\D} > 0$ on $(0, \infty)^{\kappa}_{\T}$, we have $P^{\D} - Q^{\D} > 0$ and thus, invertible on $(0, \infty)^{\kappa}_{\T}$. 
Define $B:= \left( \begin{array}{cc} a(t) & b(t) \\ b(t) & a(t) \end{array} \right)$ with $a(t)>b(t)$ for all $t \in (0, \infty)_{\T}$. Then $B$ and any real symmetric matrix of the form $Q$ will commute with $e_{B}(t,0)$, as there exists an orthogonal matrix $M = \left( \begin{array}{cc} 1 & 1 \\ -1 & 1 \end{array} \right)$ such that $M^{-1}BM$, $M^{-1}e_{B}(t,0)M$ and $M^{-1}QM$ are diagonal matrices of their respective eigenvalues. Thus, the principal axes of the associated quadric surface of $e_{B}(t,0)$ coincide with the principal axes of the associated quadric surfaces of $B$ and $Q$ (see \cite[p.7]{pipes}). 

Therefore, taking $a=u^{\D}$ and $b=0$, we obtain, for all $t \in (0, \infty)_{\T}$, 
\begin{eqnarray*}
(P-Q)^{-1} - (P^{\D} - Q^{\D})^{-1} B(t) &=& \left( \begin{array}{cc} 1/u - 1 & 1 \\ -1 & 1/v - 1 \end{array} \right)
\end{eqnarray*}
and, for any non--zero ${\bf z} \in \R^2$ with $z= \left( \begin{array}{c} x \\ y \end{array} \right)$, we have, for all $t \in (0, \infty)_{\T}$,
\begin{eqnarray*}
{\bf z}^{T} [(P-Q)^{-1} - (P^{\D} - Q^{\D})^{-1} B(t)] {\bf z} &=& \dfrac{1-u}{u} x^2 + \dfrac{1-v}{v} y^{2}\\
&<& 0. 
\end{eqnarray*}
Thus, we have $(P-Q)^{-1} < (P^{\D} - Q^{\D})^{-1} B(t)$ for all $t \in (0, \infty)_{\T}$ by Theorem \ref{propposdef}(7). This completes all conditions of Corollary \ref{exis4} and we conclude that \eqref{1}, \eqref{1i} has, at most, one positive definite symmetric solution, $X$, on $(0, \infty)_{\T}$.
\begin{flushright}
$\square$
\end{flushright}

Our next result concerns the non--multiplicity of solutions to the dynamic IVPs \eqref{2}, \eqref{1i} which Theorem \ref{exis1} or Corollary \ref{exis2} do not directly apply to. However, we employ the regressiveness of a positive definite matrix $B$ to prove the non--multiplicity of solutions to the IVP \eqref{2}, \eqref{1i}, within a domain $S \subseteq \R^{n^2}$ by constructing a modified Lipschitz condition.

\begin{theorem} \label{exis3} 
Let $S \subseteq \R^{n^2}$ and let $F:[a,b]_\T \times S \rightarrow \R^{n^2}$ be a right--Hilger--continuous function. If there exist $P, Q \in S$ with $P>Q$ on $[a,b]_{\T}$ and a positive definite matrix $B$ on $\T$ such that 
\begin{itemize}
\item[(1)] $B \in C_{rd}([a,b]_{\T};\R^{n^2})$;
\item [(2)] $e_{B}(t,a)$ commutes with $B(t)$ for all $t \in [a,b]_{\T}$ and with $P(t)$ for all $(t,P) \in [a,b]_{\T} \times S$;
\item[(3)] the inequality
\begin{equation} \label{LCOS3}
F(t,P) - F(t,Q) \le - \ominus B(t) ( P - Q )
\end{equation}
holds for all $(t,P), (t,Q) \in [a,b]]^{\kappa}_\T \times S$,
\end{itemize}
then the IVP \eqref{2}, \eqref{1i} has, at most, one solution, $X$, with $X(t) \in S$ for all $t \in [a,b]_{\T}$.
\end{theorem}

\noindent {Proof}: As before, we consider $X, Y \in S$ as two solutions of \eqref{2}, \eqref{1i} and assume $X-Y \ge 0$ on $[a,b]_\T$. Let $W:= X - Y$. We show that $W \equiv 0$ on $[a,b]_{\T}$, and so $X(t) = Y(t)$ for all $t \in [a,b]_{\T}$.\\

Since \eqref{LCOS3} holds, we have, for all $t \in [a,b]_{\T}^{\kappa}$, 
\begin{equation}
W^{\D}(t) + \ominus B(t) W^{\si}(t) =  F(t, X^{\si}(t)) - F(t, Y^{\si}(t)) + \ominus B(t)\ (X^{\si}(t) - Y^{\si}(t)) \le 0. 
\end{equation}
Note that $I+\mu(t)B(t)$ is invertible for all $t \in [a,b]_\T$. Then, by Theorem \ref{propoplus}{(2)}, the above inequality reduces to
\begin{equation}
W^{\D}(t) - [I + \mu(t)B(t)]^{-1} B(t) W^{\si}(t) \le 0, \qquad \mbox{for all} \ t \in [a,b]_{\T} \label{ineq1}
\end{equation}
Also, $e_{B}(t,a)$ and $e_{B}^{\si}(t,a)$ are positive definite and hence invertible on $[a,b]_\T$ and, thus, from Theorem \ref{propexp}{(2)}
\begin{eqnarray}
W^{\D}(t) - e_{B}(t,a)(e_{B}^{\si}(t,a))^{-1} B(t)\ W^{\si}(t) &\le& 0, \qquad \mbox{for all} \ t \in [a,b]_{\T}. \label{last}
\end{eqnarray}
By virtue of $\emph{(2)}$, $e_{B}^{-1}(t,a)$ commutes with $B(t)$ and, so, $e_{B}(t,a)$ commutes with $e_{B}^{\si}(t,a)$, for all $t \in [a,b]_{\T}$. Thus, $e_{B}^{-1}(t,a)$ commutes with $(e_{B}^{\si}(t,a))^{-1}$ for all $t \in [a,b]_{\T}$. We also see from $\emph{(2)}$ that $e_{B}^{-1}(t,a)$ commutes with $P(t)$ for all $t \in [a,b]_{\T}$. Hence, $e_{B}^{-1}(t,a)$ commutes with $P^{\si}$ and $P^{\D}$ and, thus, with $W^{\D}$ and $W^{\si}$ for all $t \in [a,b]_{\T}$. Thus, rearranging inequality \eqref{last} and using Theorem \ref{propposdef}(6) yields
\begin{eqnarray} \label{ebw}
e_{B}^{-1}(t,a) W^{\D}(t) - (e_{B}^{\si}(t,a))^{-1} B(t)\ W^{\si}(t) &\le& 0, \qquad \mbox{for all} \ t \in [a,b]^{\kappa}_{\T}.
\end{eqnarray}

Hence, using properties of Theorem \ref{main},  Theorem \ref{propposdef} and Theorem \ref{propexp} and  with \eqref{ebw}, we obtain, for all $t \in [a,b]_\T$
\begin{eqnarray*}
[e_{B}^{-1}(t,a) W(t)]^\D &=& e_{B}^{-1}(t,a) W^\D(t) + [e_{B}^{-1}(t,a)]^\D W^\si(t) \\
&\le& e_{B}^{-1}(t,a) W^\D(t) - [e_{B}^{\si}(t,a)]^{-1} B(t) W^\si(t) \\
&\le& 0.
\end{eqnarray*}
Thus $e_{B}^{-1}(t,a) W(t)$ is non--increasing for all $t \in [a,b]_{\T}$. Since $e_{B}^{-1}(t,a) > 0$ for all $t \in [a,b]_{\T}$ and $W(a) = 0$, we have $W \equiv 0$ on $[a,b]_{\T}$. This means that $X(t) = Y(t)$ for all $t \in [a,b]_{\T}$.
\begin{flushright}
$\square$
\end{flushright}

\begin{theorem} \label{exis3right} 
Let $S \subseteq \R^{n^2}$ and let $F:[a,b]_\T \times S \rightarrow \R^{n^2}$ be a right--Hilger--continuous function. If there exist $P, Q \in S$ with $P > Q$ on $[a,b]_{\T}$ and a positive definite matrix ${\cC}$ on $\T$ such that
\begin{itemize}
\item[(1)] ${\cC} \in C_{rd}([a,b]_{\T};\R^{n^2})$;
\item [(2)] $e_{{\cC}}(t,a)$ commutes with ${\cC}(t)$ for all $t \in [a,b]_{\T}$ and with $P(t)$ for all $(t,P) \in [a,b]_{\T} \times S$;
\item[(3)]  the inequality
\begin{equation} \label{LCOS3right}
F(t,P) - F(t,Q) \le ( P - Q ) (-\ominus {\cC}(t))
\end{equation}
holds for all $(t,P), (t,Q) \in [a,b]]^{\kappa}_\T \times S$,
\end{itemize}
then the IVP \eqref{2}, \eqref{1i} has, at most, one solution, $X$, with $X(t) \in S$ for all $t \in [a,b]_{\T}$.
\end{theorem}

\noindent {Proof}: The proof is similar to that of Theorem \ref{exis3} and is omitted.
\begin{flushright}
$\square$
\end{flushright}

\begin{corollary} \label{exis5} 
Let $S \subseteq \R^{n^2}$ and $F: [a,b]^{\kappa}_{\T} \times S \rightarrow \R^{n^2}$ be right--Hilger--continuous. Assume there exists a positive definite matrix $B$ on $\T$ such that conditions (1) and (2) of Theorem \ref{exis3} hold.  If $P-Q$ is positive definite and increasing on $[a,b]_{\T}$ and the inequality 
\begin{equation} \label{LCOS5}
(P-Q)^{-1} \le (P^{\D} - Q^{\D})^{-1} (- \ominus B) , \qquad \mbox{for all} \ (t,P), (t,Q) \in [a,b]]^{\kappa}_{\T} \times S
\end{equation}
holds, then the IVP \eqref{2}, \eqref{1i} has, at most, one solution $x$ with $x(t) \in S$ for all $t \in [a,b]_{\T}$.
\end{corollary}

\noindent {\bf Proof:} If \eqref{LCOS5} holds then \eqref{LCOS3} holds, by Theorem \ref{propposdef}(7). Hence, the IVP \eqref{2}, \eqref{1i} has, at most, one solution by Theorem \ref{exis3}.
\begin{flushright}
$\square$
\end{flushright}

\begin{corollary} \label{exis5right} 
Let $S \subseteq \R^{n^2}$ and $F: [a,b]^{\kappa}_{\T} \times S \rightarrow \R^{n^2}$ be right--Hilger--continuous. Assume there exists a positive definite matrix ${\cC}$ on $\T$ such that conditions (1) and (2) of Theorem \ref{exis3right} hold.  If $P-Q$ is positive definite and increasing on $[a,b]_{\T}$ and the inequality 
\begin{equation} \label{LCOS5right}
(P-Q)^{-1} \le  -\ominus {\cC} (P^{\D} - Q^{\D})^{-1}  , \qquad \mbox{for all} \ (t,P), (t,Q) \in [a,b]_{\T} \times S
\end{equation}
holds, then the IVP \eqref{2}, \eqref{1i} has, at most, one solution $x$ with $x(t) \in S$ for all $t \in [a,b]_{\T}$.
\end{corollary}

\noindent {\bf Proof:} If \eqref{LCOS5right} holds then \eqref{LCOS3right} holds, by Theorem \ref{propposdef}(6). Hence, the IVP \eqref{2}, \eqref{1i} has, at most, one solution by Theorem \ref{exis3right}.
\begin{flushright}
$\square$
\end{flushright}

We present an example of a matrix dynamic equation that has a unique solution, using Theorem \ref{exis3} and the following lemma \cite[Theorem 5.27]{BP}.

\begin{lemma} \label{linear}
Let $a,b \in \T$ with $b > a$ and $X:[a,b]_{\T} \to \R^{n^2}$. Consider the matrix initial value problem
\begin{eqnarray*}
X^{\D} &=& -V^{*}(t) X^{\si} + G(t), \qquad \ \mbox{for all} \ t \in [a,b]_{\T}; \\
X(a) &=& A,
\end{eqnarray*}
where $G$ is a rd--continuous $n^2$--matrix function on $[a,b]_{\T}$. If $V:[a,b]_{\T} \to \R^{n^2}$ is regressive then the above IVP has a unique solution $$ X(t) = e_{\ominus V^{*}}(t,a) A + \int_{a}^{t} e_{\ominus V^{*}}(t,s) G(s) \D s, \qquad \mbox{for all} \ t \in [a,b]_{\T}.$$
\end{lemma}

\begin{example} 
Let $S$ be the set of all non--singular symmetric matrices of order $n^2$. Let $K = a_{i}I$ for $1 \le i \le n$, where $a_{i} \in (0, \infty)_{\T}$.\\

\noindent Consider the IVP
\begin{eqnarray} 
X^\D = F(t, X^{\si}) &=& -K(I+2\mu(t)K)^{-1} X^{\si} + e_{{\tiny \ominus K (I + 2\mu(t)K)^{-1}}}(t,a), \label{ex2}\\ 
&& \mbox{for all} \ t \in [a, b]^{\kappa}_{\T}; \nonumber \\
X(a) &=& I. \label{ex2i}
\end{eqnarray}
We claim that \eqref{ex2}, \eqref{ex2i} has, at most, one solution, $X$, such that $X \in S$ for all $t \in [a, b]_{\T}$.
\end{example}

\noindent {\bf Proof}: We note that $K$ is a positive definite and diagonal matrix and hence $I+\mu K$ is invertible on $[a,b]_{\T}^{\kappa}$ and commutes with $K$. Moreover, $-K(I+\mu K)^{-1}$ is also diagonal and, thus, commutes with $P$.

We also note that $F$ is right--Hilger--continuous on $[a,b]_{\T} \times \R^{n^2}$, as each of its components is rd--continuous on $[a,b]_{\T}$. It follows from Theorem \ref{propoplus}(2) that, for all $t \in [a,b]_{\T}$, we have
\begin{eqnarray}
F(t,P) - F(t,Q) + \ominus 2K(P-Q) &=& [-K(I + 2\mu(t)K)^{-1}-2K(I + 2\mu(t)K)^{-1}] (P - Q)\nonumber \\
&=& -3K(I + 2\mu(t)K)^{-1}(P - Q)\\
&<& 0,
\end{eqnarray}
where we used \ref{propposdef}(6) in the last step. Therefore, \eqref{LCOS3} holds for $B=2K$. Hence, the IVP \eqref{ex2}, \eqref{ex2i} has, at most, one solution $X$ such that $X \in S$. 

Moreover, by Theorem \ref{linear}, the non--singular matrix function $$ X(t) = e_{{\tiny \ominus K(I + 2\mu(t)K)^{-1}}}(t,a) (1+ t-a)$$ uniquely solves \eqref{ex2}, \eqref{ex2i} for all $t \in [a,b]_{\T}$.

\section{Conclusions and future directions}

In this paper, we presented results identifying conditions that guarantee that if the systems \eqref{1}, \eqref{1i} and \eqref{2}, \eqref{1i} have a solution then it is unique. We did this by formulating suitable Lipschitz conditions for matrix--valued functions on time scales. The conditions will also be helpful to determine the existence and uniqueness of solutions to dynamic models of the form \eqref{1}, \eqref{1i} and \eqref{2}, \eqref{1i} and of the higher order. The results will also be helpful to establish properties of solutions for matrix--valued boundary value problems on time scales.

\section{Acknowledgements}
The author is grateful to Dr Chris Tisdell for his useful questions and comments that helped to develop the ideas in this work.


\begin{thebibliography}{12}
\markboth{Taylor \& Francis and I.T. Consultant}{Journal of Difference Equations and Applications}

\bibitem{AL}
R. P. Agarwal, and V. Lakshmikantham. {\em Uniqueness and nonuniqueness criteria for ordinary differential equations}, World Scientific Publishing Co. Inc., River Edge, NJ, 1993.



\bibitem{BS}
M. Berzig, and B. Samet. {\em Solving systems of nonlinear matrix equations involving {L}ipshitzian mappings}, Fixed Point Theory Appl. 89 (2011), 10, pp.1687-1812.

\bibitem{bhatia}
R. Bhatia. {\em Positive definite matrices}, Princeton Series in Applied Mathematics, Princeton University Press, Princeton, NJ, 2007. 

\bibitem{Birk}
G. Birkhoff, and G. Rota. {\em Ordinary differential equations}, Fourth Edition, John Wiley \& Sons Inc., New York, 1989.

\bibitem{BP}
M. Bohner, and A. Peterson. {\em Dynamic equations on time scales: An introduction with applications}, Birkh\"auser Boston Inc., Boston, MA, 2001.

\bibitem{EC}
E. A. Coddington. {\em An introduction to ordinary differential equations}, Prentice-Hall Mathematics Series, Prentice-Hall Inc., Englewood Cliffs, N.J., 1961.

\bibitem{Hart}
P. Hartman. {\em Ordinary differential equations}, John Wiley \& Sons Inc., New York, 1964.

\bibitem{Hilger}
S. Hilger. {\em Analysis on measure chains--a unified approach to continuous and discrete calculus}, Results Math., Results in Mathematics, 18 (1990),  No.1--2, pp. 18--56.

\bibitem{Hilger2}
S. Hilger. {\em Differential and difference calculus---unified!}, Proceedings of the {S}econd {W}orld {C}ongress of {N}onlinear {A}nalysts, {P}art 5 ({A}thens, 1996), Nonlinear Analysis. Theory, Methods \& Applications, 30 (1997), No. 5, pp 2683--2694.

\bibitem{Horn}
R. A. Horn, and C. R. Johnson. {\em Matrix analysis}, Cambridge University Press, Cambridge, 1990.


\bibitem{KP}
W. Kelly, and A. Peterson. {\em The theory of differential equations classical and qualitative}, Pearson Education, Inc., Upper Saddle River, NJ 07458, 2004.

\bibitem{pipes}
L. A. Pipes. {\em Matrix methods for engineering}, Prentice-Hall Inc., Englewood Cliffs, NJ, 1963.


\bibitem{KRP}
K. R. Prasad. {\em Matrix Riccati differential equations on time scales}, Comm. Appl. Nonlinear Anal., 8 (2001), No 4, pp 63--75.


\bibitem{TZ1}
C. C. Tisdell, and A. H. Zaidi. {\em Successive approximations to solutions of dynamic equations on time scales}, Comm. Appl. Nonlinear Anal., 16 (2009), No. 1, pp. 61--87.

\bibitem{TZ}
C. C. Tisdell, and A. Zaidi. {\em Basic qualitative and quantitative results for solutions to nonlinear, dynamic equations on time scales with an application to economic modelling}, Nonlinear Analysis. Theory, Methods \& Applications, 68 (2008), No. 11, pp. 3504--3524.

\bibitem{thesis}
A. Zaidi. {\em Existence and Uniqueness of solutions to nonlinear first order dynamic equations on time scales}, Ph.D. thesis, University of New South Wales, 2009.

\bibitem{zaidi}
A. H. Zaidi. {\em Existence of solutions and convergence results for dynamic initial value problems using lower and upper solutions}, Electronic Journal of Differential Equations, 2009(2009), No. 161, pp. 1--13.

\end{thebibliography}
\end{document}